\documentclass[12pt]{amsart} 
\usepackage{color, amsfonts, amssymb,geometry}
\usepackage{latexsym}
\usepackage{graphicx}
\usepackage{amssymb}

\newcommand{\IS}{\mathbb{S}}
\newcommand{\R}{\mathbb{R}}

\newcommand{\Z}{\mathbb{Z}}
\newcommand{\C}{\mathbb{C}}

\newcommand{\IH}{\mathbb{H}}

\newcommand{\rs}{\mbox{$\widehat{\C}$}}

\def\DDD{{\mathcal D}}
\newcommand{\ARconfdim}{\hbox{ARconfdim}}
\newcommand{\WF}{\mathfrak{P}}
\newcommand{\vp}{\varphi}

\newcommand{\pf}{\noindent {\bf Proof: }}
\newcommand{\bdry}{\partial}                     
\newcommand{\intersect}{\cap}                    
\newcommand{\union}{\cup}                        
\newcommand{\cl}{\overline}                      

\newcommand{\gap}{\vspace{5pt}}                 

\newcommand{\be}{\begin{enumerate}}
\newcommand{\eb}{\end{enumerate}}

\newtheorem{prop}{Proposition}[section]

\newtheorem{lemma}{Lemma}[section]
\newtheorem{cor}{Corollary}[section]

\newtheorem{mainthm}{Theorem}

\title{Quasisymmetrically inequivalent hyperbolic Julia sets}
\author{Peter Ha\"{\i}ssinsky}
\address{LATP/CMI\\ Universit\'e de Provence\\ 39, rue Fr\'ed\'eric Joliot-Curie\\13453 Marseille cedex 13\\France}                        
\email{phaissin@cmi.univ-mrs.fr}
\author{Kevin M. Pilgrim}
\address{Dept. of Mathematics\\ Indiana University\\Bloomington\\ IN 47405\\  USA}
\email{pilgrim@indiana.edu}
\date{\today}                                           

\begin{document}
\maketitle

\begin{abstract}
We give explicit examples of pairs of Julia sets of hyperbolic rational maps
which are homeomorphic but not quasisymmetrically homeomorphic. 
\end{abstract}

\gap
\gap

Quasiconformal geometry is concerned with properties of metric spaces that are preserved under 
quasisymmetric homeomorphisms. Recall that a homeomorphism $h: X \to Y$ between metric spaces is {\em quasisymmetric} 
if there exists a distortion control function $\eta: [0,\infty) \to [0,\infty)$ which is a homeomorphism and which satisfies 
$|h(x)-h(a)|/|h(x)-h(b)| \leq \eta(|x-a|/|x-b|)$ for every triple of distinct points $x, a, b \in X$. 
We shall say that $X$ and $Y$ are quasisymmetrically equivalent if there exists such a homeomorphism.   

\gap

A basic ---even if still widely open--- question is to determine whether two given spaces belong to the same quasisymmetry
class, once it is known that they are homeomorphic and share the same qualitative geometric properties. This question
arises also in the classification of hyperbolic spaces and word hyperbolic groups in the sense of Gromov \cite{bourdon:pajot:survey,
kleiner:icm2006,ph:bbki}.
Besides spaces modelled on manifolds, very few examples are understood; see nonetheless \cite{bourdon:gafa} for examples
of inequivalent spaces modelled on the universal Menger curve.
Here, we focus our attention on compact metric spaces that arise as Julia sets of rational maps.  
A rational map is {\em hyperbolic} if the closure of the set of forward orbits of all its critical points 
does not meet its Julia set.  
We address the question of whether the geometry of the Julia set of a hyperbolic rational map
is determined by its topology. More precisely, {\it given two hyperbolic rational maps $f$ and $g$ with
homeomorphic Julia sets $J_f$ and $J_g$, does there exist a quasisymmetric homeomorphism $h:J_f\to J_g$?}

\gap

Hyperbolic Julia sets serve our purposes for several reasons. First, 
it rules out elementary local obstructions. For instance, the Julia set of $f(z)=z^2$ is the Euclidean unit circle $\IS^1$, 
while that of  $g(z)=z^2+1/4$ is a Jordan curve with a cusp at the unique fixed-point, so they are not quasisymmetrically
equivalent.
Second, if $f$ is hyperbolic, it is locally invertible near $J_f$, and the inverse branches are uniformly contracting; 
the Koebe distortion theorem then implies that $J_f$ satisfies a strong quasi-self-similarity property.   Among such maps, in some cases, this implies that homeomorphic Julia sets are quasisymmetrically homeomorphic.  

\be
\item If the Julia set of a hyperbolic rational map is a Jordan curve, then it is quasisymmetrically equivalent
to the unit circle \cite{DS1}.
\item Let $C \subset \R$ denote the usual middle-thirds Cantor set.  Recall that any compact, totally disconnected metric space 
without isolated points is homeomorphic to $C$; see e.g. \cite[Thm. 2.97]{HockingYoung}. 
If the Julia set of a hyperbolic rational map is homeomorphic to $C$, then, by a theorem of 
David and Semmes \cite[Prop.\,15.11]{david:semmes:dreams} 
they are quasisymmetrically equivalent.  
\item If $f$ and $g$ are hyperbolic and their Julia sets are homeomorphic by the restriction of
a global conjugacy, then they are also quasisymmetrically equivalent \cite[Thm 2.9]{ctm:ds:qciii}.
\eb

So one must look to more complicated Julia sets 
for potential examples of nonquasisymmetrically equivalent Julia sets.  

We will show 

\begin{mainthm}
\label{thm:ctimesq}
Let $f(z)=z^2+10^{-9}/z^3$ and $g(z)=z^2+10^{-20}/z^4$.  Then $J_f, J_g$ are each homeomorphic to $C \times S^1$, 
but they are not quasisymmetrically homeomorphic.  
\end{mainthm}

Recall that a metric space $X$ equipped with a Radon measure $\mu$ is {\em Ahlfors regular of dimension $Q$} 
if the measure of a ball satisfies $\mu(B(x,r)) \asymp r^Q$; one has then that $X$ has locally finite 
Hausdorff measure in its Hausdorff dimension, $Q$.  Its {\em Ahlfors-regular conformal dimension} $\ARconfdim(X)$ 
is the infimum of the Hausdorff dimensions of all Ahlfors-regular metric spaces quasisymmetrically equivalent to $X$.  
Since the Julia set of any hyperbolic rational map is quasi-self-similar, it follows that it is Ahlfors regular
and porous, hence has Hausdorff dimension strictly less than $2$ \cite[Thm 4 and Cor.]{DS2}.  
So if $f$ is hyperbolic then $\ARconfdim(J_f)<2$.  
We prove Theorem \ref{thm:ctimesq} by showing  $\ARconfdim(J_f) \neq \ARconfdim(J_g)$.  

\gap

The arguments we use to prove Theorem \ref{thm:ctimesq} will generalize to yield

\begin{mainthm}
\label{thm:ctimesq_sequence}
There exist hyperbolic rational maps each of whose Julia sets is homeomorphic to $C \times S^1$ and 
whose Ahlfors-regular conformal dimensions are arbitrarily close to $2$.  
\end{mainthm}

It follows that there exists an infinite sequence of hyperbolic rational maps whose Julia sets are homeomorphic to 
$C \times S^1$ but which are pairwise quasisymmetrically inequivalent.  

\gap

Our method of proof of Theorem \ref{thm:ctimesq_sequence} requires that the degrees become arbitrarily large.  
It is tempting to look for such a sequence of examples among maps of fixed degree.  
This may be difficult: as is shown by Carrasco \cite{carrasco:phd}, the Ahlfors-regular conformal dimension of any hyperbolic polynomial 
with connected Julia set is equal to $1$.  

If connected, the Julia sets of hyperbolic polynomials have many cut-points.  At the opposite extreme, 
recall that a Sierpi\'nski carpet may be defined as 
a one-dimensional, connected, locally connected compact subset of the sphere such that the components of its complement
are Jordan domains with pairwise disjoint closures; any two such spaces are homeomorphic \cite{whyburn:sierpinski}. 
Sierpi\'nski carpets are 1-dimensional analogs of Cantor sets. They also play an important role in Complex Dynamics and Hyperbolic Geometry \cite{ctm:classification,bonk:icm:qcgeom}.
Sierpi\'nski carpets which arise from hyperbolic groups and hyperbolic rational maps also share the same qualitative properties:
their peripheral circles are uniform quasicircles and are uniformly separated; Bonk also proved that  any such carpet is quasisymmetrically equivalent to one where the complementary domains are round disks in $\rs$.   Nonetheless, using similar methods, we will show 

\begin{mainthm}
\label{thm:carpets}
There exist hyperbolic rational maps with Sierpi\'nski carpet Julia sets whose Ahlfors-regular conformal dimensions 
are arbitrarily close to $2$.
\end{mainthm}

To our knowledge, an analogous result for conformal dimensions of limit sets of convex compact Kleinian groups is not yet known.  

On the one hand, it is perhaps not surprising that there are a plethora of quasisymmetrically distinct such Julia sets: any quasisymmetric map 
between round convex compact Kleinian group carpets is the restriction of a M\"obius transformation \cite[Thm. 1.1]{bonk:kleiner:merenkov:schottky}.  
Also,  any quasisymmetric automorphism of the standard square ``middle ninths'' carpet is the restriction of a Euclidean isometry 
\cite[Thm. 8.1]{bonk:icm:qcgeom}.  On the other hand, the proofs of these results are rather involved.

The proofs of our results rely on the computation of the Ahlfors-regular conformal dimension of 
certain metric spaces homeomorphic to $C\times S^1$, following the seminal work of Pansu, cf. \cite[Prop.\,3.7]{ph:bbki}.  
We will also make frequent use of the fact that on the Euclidean $2$-sphere, an orientation-preserving self-homeomorphism is quasiconformal if and only if it is quasisymmetric; 
see  \cite[Thm. 11.14]{heinonen:analysis}. We denote by $\IS^2$ denote the round Euclidean $2$-sphere.

The special case needed for the present purpose is summarized in \S\,1. The proofs of the theorems appear in \S\S\,2 and 3.  

\section{Annulus maps}

Let $I=[0,1]$ and let $\iota: I \to I$ be the involution given by $\iota(x)=1-x$.  Identify $\IS^1$ with $\R/\Z$.  
We give $I \times \IS^1$ the product orientation.  Fix an even integer $m\geq 2$ and let $\DDD:=(d_0, \ldots, d_{m-1})$ 
be a sequence of positive integers such that $\sum_{i=0}^{m-1} \frac{1}{d_i} < 1$.  
Then there exist real numbers $a_i, b_i, i=0, \ldots, m-1$ such that for each $i$, $|b_i-a_i|=\frac{1}{d_i}$ and 
\[ 0 < a_0 < b_0 < a_1 < b_1 < \ldots < a_{m-1} < b_{m-1} < 1.\]
Fix such a choice $a_0, b_0, \ldots, a_{m-1}, b_{m-1}$.  
For each $i$, let $J_i=[a_i, b_i]$, and let $g_i: I \to J_i$ be the unique affine homeomorphism 
which is orientation-preserving if $i$ is even and is orientation-reversing if $i$ is odd.  
This iterated function system on the line has a unique attractor $C(\DDD)$ and its Hausdorff dimension, 
by the pressure formula \cite[Thm. 5.3]{falconer:techniques}, is the unique real number $\lambda=\lambda(\DDD)$ satisfying 
\[ \sum_{i=0}^{m-1}\frac{1}{d_i^\lambda} = 1.\]

Let 
\[ \tilde{F}: \left(\sqcup_{i=0}^{m-1}J_i\right)\times \IS^1 \to I \times \IS^1 =: A. \]
be the map whose restriction to the annulus $A_i:=J_i \times \IS^1$ is given by 
\[ \tilde{F}|_{A_i}(x, t) = (g_i^{-1}(x) , (-1)^i d_i \cdot t \bmod 1).\]

That is:  $\tilde{F}|_{A_i}$ is an orientation-preserving covering map of degree $d_i$ which is a Euclidean homothety with factor 
$d_i$ and which preserves or reverses the linear orientation on the interval factors according to whether $i$ is even or, 
respectively, is odd.  

The invariant set associated to $\tilde{F}$ is 
\[ X(\DDD) := C(\DDD) \times \IS^1 = \bigcap_{n\geq 0} \tilde{F}^{-n}(A).\]
From \cite[\S 3]{kmp:ph:examples}, we have 

\begin{prop}
\label{prop:compute_confdim}
The Ahlfors-regular conformal dimension of $X(\DDD)$ is equal to $1+\lambda(\DDD)$.  
\end{prop}

This statement is a particular case of a criterion originally due to Pansu \cite[Prop. 3.7]{ph:bbki}; see also
Tyson's theorem \cite[Thm 15.10]{heinonen:analysis}.

\section{Proofs of Theorems \ref{thm:ctimesq} and \ref{thm:ctimesq_sequence}}

Let $\DDD$ be a sequence of positive integers defining a family of annulus maps $\tilde{F}$ as in the previous section, 
and put $X=X(\DDD)$.  

\begin{prop} 
\label{prop:extend}
There is a smooth embedding $A \hookrightarrow \IS^2$ such that (upon identifying $A$ with its image) the map 
$\tilde{F}: \sqcup_i A_i \to A$ extends to a smooth map  $F: \IS^2 \to \IS^2$ whose iterates are uniformly quasiregular.  
There is a quasiconformal (equivalently, a quasisymmetric) homeomorphism $h: \IS^2 \to \rs$ such that 
$h\circ F \circ h^{-1}$ is a hyperbolic rational map $f$, and $h(X)=J_f$.  
\end{prop}

\pf  The existence of the extension $F$ is a straightforward application of quasiconformal surgery; 
we merely sketch the ideas and refer to \cite{kmp:tan:surgery} for details; see also the forthcoming text \cite{branner:surgery} 
devoted to this topic.  The next two paragraphs outline this construction.

The linear ordering on the interval $I$ gives rise to a linear ordering on the set of $2m$ boundary components of the set 
of annuli $A_0, \ldots, A_{m-1}$.  
We may regard $A$ as a subset of a smooth metric sphere $S^2$ conformally equivalent to $\IS^2$.  
For $i=1, \ldots, m-1$ let $C_i$ be the annulus between $A_{i-1}$ and $A_i$.  
Let $D_0, D_1$ be the disks bounded by the least, respectively greatest, boundary of $A$, so that the interiors of 
$D_0, A, D_1$ are disjoint.   
Let $D_0'$ be the disk bounded by the least component of $A_0$ and $D_1'$ be the disk bounded by the greatest component of 
$A_{m-1}$.  

We now extend $\tilde{F}$ as follows.  See Figure 1.
\begin{figure}
\label{fig:caric}
\begin{center}
\includegraphics[width=5in]{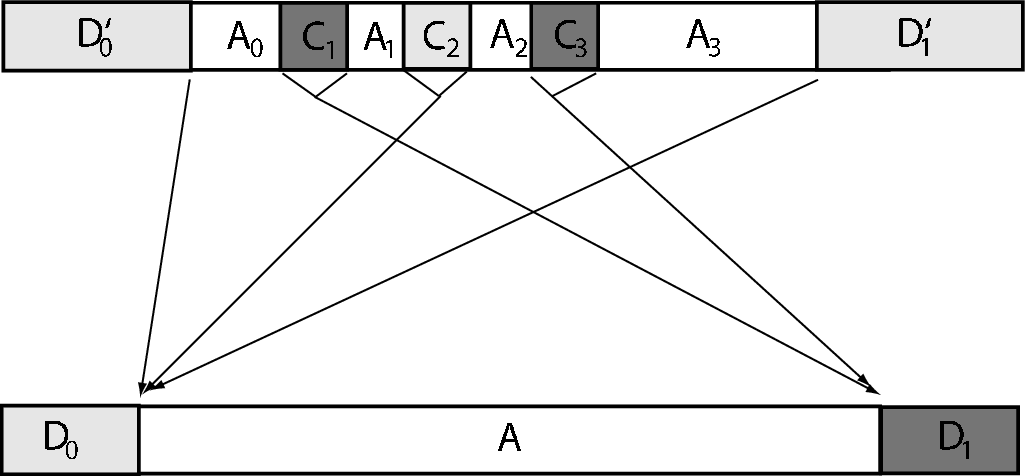}
\caption{Caricature of the extended mapping, $F$. }
\end{center}
\end{figure}
Send $D_0' $ to $D_0$ by a proper map of degree $d_0$ ramified over a single point $x$, so that in suitable holomorphic coordinates 
it is equivalent to $z \mapsto z^{d_0}$ acting near the origin; thus $D_0 \subset D_0'$ is mapped inside itself.  
Similarly, send  $D_1' $ to $D_0$ by a proper map of degree $d_{m-1}$ ramified only over $x$, so that in suitable holomorphic 
coordinates it is equivalent to $z \mapsto 1/z^{d_{m-1}}$ acting near infinity; thus $D_1 \subset D_1'$ is mapped into $D_0$.  
To extend over the annulus $C_i$ between $A_{i-1}$ and $A_{i}$, note that both boundary components of $C_i$ map either to the least, 
or to the greatest, component of $\bdry A$.  It is easy to see that there is a smooth proper degree $d_{i-1}+d_{i}+1$ 
branched covering of $C_i$ to the corresponding disk $D_0$ (if $i$ is even) or $D_1$ (if $i$ is odd).  
This completes the definition of the extension $F$.

It is easy to arrange that $F$ is smooth, hence quasiregular.  We may further arrange so that the locus where $F$ is not conformal 
is contained in a small neighborhood of $C_1 \union \ldots \union {C_{m-1}}$.  This locus is nonrecurrent, so the iterates of 
$F$ are uniformly quasiregular.  
By a theorem of Sullivan \cite[Thm 9]{DS2}, $F$ is conjugate via a quasiconformal homeomorphism $h: S^2 \to \rs$  
to a rational map $f$.  By construction, every point not in $h(X)$ converges under $f$ to a superattracting fixed-point $h(x)$ in  
the disk $h(D_0)$, so $f$ is hyperbolic and $h(X)=J_f$.   \qed

We now establish a converse.

\begin{prop}
\label{prop:restrict}
Suppose $f: \rs \to \rs$ is a rational map for which there exists a closed annulus $A$ and essential pairwise disjoint subannuli 
$A_0, A_1, \ldots, A_{m-1}$, $m$ even, contained in the interior of $A$ such that (with respect to a linear ordering induced by $A$) 
$A_0 < A_1 < \ldots A_{m-1}$.  Let $D_0, D_1$ be the disks bounded by the least (respectively, greatest) boundary component of $A$. 
Further, suppose that for each $i=0, \ldots, m-1$, $f|_{A_i}: A_i \to A$ is a proper covering map of degree $d_i$, 
with $f$ mapping the greatest component of $A_i$ and the least component of $A_{i+1}$ to the boundary of $D_1$ if $i$ is even, 
and to the boundary of $D_0$ if $i$ is odd.   
Put $\DDD = (d_0, d_1, \ldots, d_{m-1})$.  Let $\tilde{f}=f|_{\sqcup_{i=0}^{m-1}A_i}$ and put 
$Y=\intersect_{n\geq 0}\tilde{f}^{-n}(A)$.  Then $Y \subset J_f$, $\tilde{f}(Y)=Y=\tilde{f}^{-1}(Y)$,  
and there is a quasisymmetric homeomorphism $h: Y \to X$ conjugating $\tilde{f}|_Y: Y \to Y$ to $\tilde{F}|_X: X \to X$ 
where $\tilde{F}$ is the family of annulus maps defined by the data $\DDD$.  
\end{prop}

\pf The conformal dynamical systems of annulus maps defined by $\tilde{f}$ and by $\tilde{F}$ are combinatorially equivalent 
in the sense of McMullen \cite[Appendix A]{ctm:siegel}, so by \cite[Thm. A.1]{ctm:siegel} there exists a quasiconformal 
(hence quasisymmetric) conjugacy $\tilde{h}$ from $\tilde{f}$ to $\tilde{F}$; we set $h=\tilde{h}|_Y$.   
\qed

Combined with Proposition \ref{prop:compute_confdim}, this yields:
\begin{cor}
\label{cor:lower_bound}
Under the assumptions of Proposition \ref{prop:restrict},
$\ARconfdim(J_f) \geq 1+\lambda(\DDD)$, with equality if $Y=J_f$.  
\end{cor}

\noindent{\bf Proof of Theorem \ref{thm:ctimesq}.}  For $\epsilon \in \C$ let $f_\epsilon(z)=z^2+\epsilon/z^3$.   
McMullen \cite[\S 7]{ctm:automorphisms} shows that for $|\epsilon|$ sufficiently small the map $f_\epsilon$ restricts 
to a family of annulus maps  with the combinatorics determined by the data $\DDD=(2,3)$ and with Julia set homeomorphic 
to the repellor $X_{(2,3)}$ determined by $\DDD$; it is easy to see that $\epsilon=10^{-9}$ will do.  

Exactly the same arguments applied to the family $g_\epsilon (z)=z^2+\epsilon/z^4$ show that  if $|\epsilon|$ 
is sufficiently small, the family $g_\epsilon$ restricts to a family of annulus maps with the combinatorics determined 
by $\DDD=(2,4)$ and whose Julia set is homeomorphic to the corresponding repellor $X_{2,4}$.   
It is easy to see that $\epsilon=10^{-20}$ will do; one may take $A=\{10^{-6} < |z| < 10^{10}\}$.  
By Corollary \ref{cor:lower_bound} and Proposition \ref{prop:compute_confdim}, the Ahlfors-regular conformal dimensions 
$1+\lambda_f, 1+\lambda_g$ of $J_f, J_g$ satisfy the respective equations $2^{-\lambda_f}+3^{-\lambda_f} = 1$, 
$2^{-\lambda_g} + 4^{-\lambda_g}=1$ and are therefore unequal.  
Since the Ahlfors-regular conformal dimension is a quasisymmetry invariant, the proof is complete.\qed
\gap

\noindent{\bf Proof of Theorem \ref{thm:ctimesq_sequence}.}  For an even integer $n\geq 4$, let 
$\DDD_n = (d_0, d_1, \ldots, d_{n-1})$ 
where $d_0=\frac{1}{n+1}$ and $d_i=\frac{1}{n}$ for $i=1, \ldots, n-1$.  
Let $f_n$ be the rational map given by Proposition \ref{prop:extend}.  By Corollary \ref{cor:lower_bound} 
$\ARconfdim(J_{f_n})$ is $1$ plus the unique positive root 
$\lambda_n$ of the equation 
\[ (n+1)^{-\lambda} + (n-1)n^{-\lambda}=1.\]
The left-hand side is larger than $1$ when $\lambda=\frac{\log(n-1)}{\log(n)}$, so 
$\lambda_n > \frac{\log(n-1)}{\log n}$ and thus $\lambda_n \to 1$ as $n \to \infty$. 
Hence $\ARconfdim(J_{f_n}) \to 2$ as $n \to \infty$.  \qed

\section{Proof of Theorem \ref{thm:carpets}}

Fix an even integer $n \geq 2$.  For each such $n$, we will build a rational function $f_n: \rs \to \rs$ 
with the following properties:  (1) its  Julia set is homeomorphic to the Sierpi\'nski carpet, and 
(2) there exists an annulus $A \subset \rs$, and parallel pairwise disjoint essential subannuli 
$A_0, \ldots, A_{n-1}$ 
such that for each $i=0, \ldots, n-1$, the restriction $f|_{A_i}: A_i \to A$ is a proper holomorphic covering of 
degree $(n+4)$,  just as in the previous section.  
Theorem \ref{thm:carpets} will then follow immediately from Corollary \ref{cor:lower_bound} with 
$\DDD=(\; \underbrace{n+4, \ldots, n+4}_{n}\; )$.  

We will first build $f_n$ as a function from one Riemann sphere to another, and then re-identify domain and range.  
We are grateful to Daniel Meyer for suggesting this construction which is more explicit than our original one.   

\begin{figure}
\begin{center}
\includegraphics[width=5in]{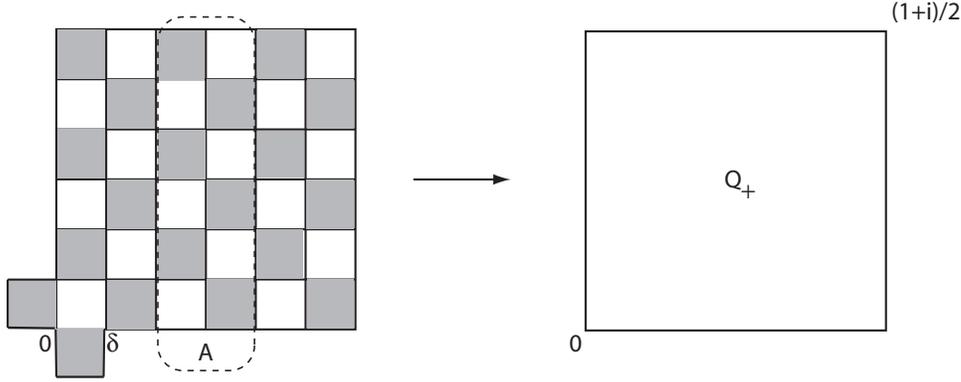}
\caption{The rational map $f_n$ when $n=2$.  The domain and codomain are the doubles of the two polygons across their boundaries.  
Note the conformal symmetries. }
\end{center}
\end{figure}

We shall suppress the dependence on $n$ in our construction.
For $z \in \rs$ set $j(z)=\bar{z}$, and let us consider the unique Weierstrass function $\WF:\C\to\rs$ which is $\Z[i]$-periodic, 
which maps $0$, $1/2$, $(1+i)/2$ and $i/2$ to $\infty$, $(-1)$, $0$ and $1$ respectively.
We may thus consider the Riemann sphere $\rs$ as the quotient Euclidean rectangle $[0,1/2]\times [(-1/2),1/2]$
upon identifying boundary-points via the map $j$. We view $\rs$ as the union of two Euclidean
squares: the ``white square'' $[0,1/2]\times [0,1/2]$ and the ``black square''
$[0,1/2]\times [(-1/2),0]$. The map $\WF$ maps the white square to the upper half-plane $\IH_+$
and the black square to the lower half-plane $\IH_-$.

To define the codomain, put $\delta = \frac{1}{2(n+4)}$, and let 
$$Q_+= [0,1/2]^2 \cup ([-\delta,0]\times [0, \delta]) \cup ([0, \delta]\times [-\delta,0])$$
and $Q_-= j(Q_+)$. Both polygons $Q_+$ and $Q_-$ are tiled by $(n+4)^2+2$ squares of size $\delta$.
Let  $\Sigma$ be the sphere obtained from the disjoint union $Q_+ \sqcup Q_-$ by gluing their boundaries via the map $j$.  Then $\Sigma$ inherits a conformal structure from that of $Q_\pm$: away from the corners this is clear; by the removable singularities theorem, 
this conformal structure extends over the corners.   Note that the map $j$ 
gives an anticonformal involution of $\Sigma$ which we denote again by $j$.

Define $F_+:Q_+\to \rs$ by $F_+(z)=\WF((n+4)z))$ and $F_-:Q_-\to\rs$ by $F_-=j\circ F_+\circ j$.
This defines a holomorphic map $F:\Sigma\to \rs$ of degree $(n+4)^2+2$. 
Considering the tiling of $\Sigma$ by the squares of size $\delta$, we may color them into
white and black in such a way that a white square is mapped under to $F$ to $\IH_+$ and
a black one to $\IH_-$, see Figure 2. The critical points of $F$ occur at where four or more squares meet.  
By construction, the image of every critical point is 
one of the points $-1, 0, 1, \infty$.

By the Riemann mapping theorem, there exists a unique conformal map $\vp_+:Q_+\to\overline{\IH_+}$
such that $\vp_+(0,1/2,(1+i)/2)=(\infty, (-1),0)$. Note that the map $x+iy\mapsto y+ix$ defines
an anticonformal involution of $Q_+$ which fixes $0$ and $(1+i)/2$ and interchanges
$1/2$ and $i/2$: this forces $\vp_+(i/2)=1$. Set $\vp_-:Q_-\to\overline{\IH_-}$
by $\vp_-=j\circ\vp_+\circ j$. Both maps patch together to form a conformal map $\vp:\Sigma\to\rs$.
Let us finally set $$f= F\circ \vp^{-1}:\rs\to\rs\,.$$
Every critical point of $f$ is first mapped to  $-1, 0, 1, \infty$, and every
point of this set maps to $\infty$ under $f$, 
which is therefore a fixed critical point at which $f$ has local degree $3$.  
Hence $f$ is a critically finite hyperbolic rational  map.

Let $A_+'= [2\delta, 1/2 - 2\delta]\times [0,1/2]\subset Q_+$ and $A_-'=j(A_+)\subset Q_-$.
Their union defines an annulus $A$ of $\Sigma$, and we let $A=\vp(A')$. 

The preimage $F^{-1}(A)$ consists of $(n+5)$ disjoint annuli, each compactly contained
in vertical strips of width $\delta$ tiled by squares. Among them, there are $n$ subannuli
$A_0',\ldots, A_{n-1}'$ compactly contained in $A'$, each map under $F$ by degree $n+4$.
Let $A_j=\vp(A_j')$: then $A_j$ is compactly contained in $A$ and $f:A_j\to A$ has degree $n+4$.
By Corollary \ref{cor:lower_bound}, $\ARconfdim(J_{f}) \geq 1+\lambda$, where $\lambda=\lambda_n$ 
is the unique positive 
root 
of the equation 
\[ n(n+4)^{-\lambda}=1.\]
As $n \to \infty$, clearly $\lambda_n \to 1$ and so $\ARconfdim(J_{f_n})\to 2$.

It remains only to show that $J_{f}$ is a Sierpi\'nski carpet.  
We imitate the arguments of Milnor and Tan Lei given in \cite[Appendix]{milnor:quadratic}.  
They first show the following:

\begin{lemma}
\label{lemma:jdomain}
Let $f$ be a hyperbolic rational map and $z$ a fixed-point at which the local degree of $f$ equals $k \geq 2$.  
Suppose $W$ is the immediate basin of attraction of $z$.  
Suppose there exist domains $U, V$ each homeomorphic to the disk such that 
$\cl{\Omega} \subset U \subset \cl{U} \subset V$ and $f|_U: U \to V$ is proper and of degree $k$.  
Then $\bdry \Omega$ is a Jordan curve.
\end{lemma} 

Note that the conformal isomorphism $\varphi: \Sigma \to \rs$ sends the union of the top and right-hand edges of the square to $[-1,1]$ and sends the point in Figure 2 labelled $0$ to infinity.  
Let $V=\rs\setminus [-1,1]$.
The map $f$ has a unique periodic Fatou component $W$ ---the immediate basin of $\infty$--- 
and clearly $W \subset V$.  
The domain $V$ is simply connected and contains exactly one critical value of $f$, namely, 
the point $\infty$.  
It follows that there is exactly one component $U$ of $f^{-1}(V)$ containing $\infty$, and 
$\cl{W} \subset U \subset \cl{U} \subset V$ and $f|_U: U \to V$ is proper and of degree $3$.  
By Lemma \ref{lemma:jdomain}, $\bdry W$ is a Jordan curve.  
The remaining arguments needed are identitical to those given in {\em op. cit.}:  
since $f$ is hyperbolic and critically finite, the Julia set is one-dimensional, connected and locally connected, and there are no critical points in the Julia set.  
It follows that every Fatou component is a Jordan domain, and that the closures of the Fatou components 
are pairwise disjoint. Therefore $J_{f}$ is homeomorphic to the Sierpi\'nski carpet  \cite{whyburn:sierpinski}.
\qed

%
%

\end{document}